\documentclass[12pt]{article}
\usepackage[all]{xy}
\usepackage{amsthm,amsfonts}
\usepackage{graphicx}
\xyoption{all}
\newtheorem{teo}{Theorem}[section]
\newtheorem{lema}[teo]{Lemma}
\newtheorem{prop}[teo]{Proposition}

\newtheorem{ex}[teo]{Example}
\newtheorem{cor}[teo]{Corollary}

\newtheorem{rmk}[teo]{Remark}
\newcommand{\sex}{\mathcal{ M}_6^b}
\newcommand{\fsex}{\overline{\mathcal{ M}}_6^b}

\newcommand{\g}{\overline{\mathcal{ H}}_g}

\def\squareforqed{\hbox{\rlap{$\sqcap$}$\sqcup$}}
\def\QED{\ifmmode\squareforqed\else{\unskip\nobreak\hfil
        \penalty50\hskip1em\null\nobreak\hfil\squareforqed
        \parfillskip=0pt\finalhyphendemerits=0\endgraf}\fi
                \medskip}
\def\p#1{\mbox{{\rm I$\!$P}$_{#1}$}}
\def\P#1{\mathbb{P}_#1}
\def\PD#1{\mathbb{P}_#1^{*}}

\def\complex{{\rm C}\hskip-6pt\raise3pt\hbox{$_{^|}$}~}
\begin{document}

 \hoffset = -1truecm
\voffset = -2truecm

\title{\textbf{Curves of genus 2 and Desargues configurations}}
  \author{
\textbf{ D. Avritzer\thanks{Both authors would like
to thank the GMD (Germany) and CNPq (Brasil)
for support during the preparation of this
paper}}
\and \textbf{  H. Lange}$^{*}$}
\date{}
\maketitle

\section{Introduction}
A Desargues configuration is the configuration of 10 points and 10 lines
of the classical theorem of Desargues in the complex projective plane. For a
precise definition see Section 2. The greek mathematician C. Stephanos
showed in 1883 (see~\cite{kn:ST}) that one can associate to every Desargues
configuration a curve of genus 2 in a canonical way. Moreover he proved
that the induced map from the moduli space of Desargues configurations to the
moduli space of curves of genus 2 is birational. Our main motivation for
writing this paper was to understand the last result. Stephanous needed about
a hundred pages of classical invariant theory to prove it. We apply 
instead a simple argument of Schubert calculus to prove a slightly more
precise version of his result.

\bigskip

Let $M_{D}$ denote the (coarse) moduli space of Desargues configurations. 
It is a three-dimensional quasiprojective variety. On the other hand, let
$\sex$ denote the moduli space of stable binary sextics. There is a canonical
isomorphism $\g \setminus \Delta_1 \cong \sex$ (see~\cite{kn:pqh}).
Here $\g$ denotes the moduli space of stable curves of genus 2 in the sense of
Deligne-Mumford and $\Delta_1$ the boundary divisor parametrizing 2 curves of
genus 1 intersecting transversely in one point. So instead of curves of genus
2, we may speak of binary sextics. 
The main result of the paper consists of the following 2 statements:\newline
(1) {\it There is a canonical injective birational morphism}
\[ \Phi :M_D \hookrightarrow \sex .\]
(2) {\it We determine a hypersurface $H \subset \sex$ such that }
\[ \sex \setminus H  \subseteq \Im m(\Phi) \subseteq
\sex, 
\]
where both the inclusions above are strict.

\bigskip

A Desargues configuration is called special if one of its lines contains
4 configuration points. We show: \newline
(3) {\it A Desargues configuration $D$ is special if and only if the binary
sextic $\Phi(D)$ admits a double point, i.e. the corresponding curve of genus
2 is not smooth.}

\bigskip

The Moduli space $M_{D}$ admits a natural compactification $\overline{M}_D.$
We also study the configurations corresponding to boundary points 
$\overline{M}_D\setminus M_D,$ which we call {\it degenerate Desargues
configurations.} In fact, there are degenerate Desargues 
configurations of the first, second and third kind (see Section 2). 
As a rational map of normal projective varieties
$\Phi:\overline{M}_D \cdots \longrightarrow \fsex$ is defined in 
codimension 1. Hence it extends to a morphism on an open set of the divisor
$\overline{M}_D\setminus M_D.$ We show, however: \newline
(4) {\it The geometric interpretation of the morphism $\Phi: \overline{M}_D
\cdots \longrightarrow \fsex$ does not extend to an open set $U:$
 $M_D \subset U \subset \overline{M}_D,$ where the first inclusion is strict.
 In fact,
the binary sextic associated to a degenerate Desargues configuration
of the first kind in an analogous way is not semistable.}

\bigskip

In Section 2, we construct the moduli spaces $M_D$ and its compactification
$\overline{M}_D.$ In Section 3, we give the definition of the map 
$\Phi: M_D \longrightarrow  \sex.$ Sections 4, 5, 6 and 7 contain the proofs
of statements (1), (2), (3) and (4) respectively. 

\bigskip

The second author would like to thank W. Barth with whom he discussed the
subject, already 15 years ago.

\section{The 
moduli space of Desargues configurations}
Let  $\P{2}$ denote the projective plane over the field of complex numbers.
The classical theorem of Desargues says:
If the lines joining corresponding vertices of two triangles $A_1,B_1,C_1$
and $A_2,B_2,C_2$ in $\P{2}$ meet in a point $A$ then the intersections of
corresponding sides lie on a line $a$ and conversely (see Figure 1).

\begin{figure}
\hspace{3.5cm}
\caption{Desargues' Theorem}
\end{figure}	

The triangles are then said to be {\it in perspective,} $A$ is called
the {\it center of perspective} and $a$ the {\it axis of perspective.}
The configuration consisting of the 10 points $A,A_i,B_i,C_i,(i=1,2,3)$
and ten lines, namely the 6 sides of the triangles, the 3 lines joining $A$  
to the vertices of the triangle and the axis of perspective, is called a 
{\it Desargues configuration.} It is a $10_3-$ configuration meaning
that each of the 10 lines contains 3 of the 10 points  and through each of the
10 points  there pass  3 of the 10 lines.

It may happen that one of the vertices of one triangle lies 
on the opposite side of the other triangle, in which case  Desargues Theorem
is still valid, but one line contains now 4 of the 10 points. Contrary to some
authors, we consider this configuration also as a Desargues configuration and
call it a {\it special Desargues configuration.}

\bigskip

As long ago as 1846, Cayley remarked (see \cite{kn:cay}, p.318) that the 10
lines and 10 planes determined by 5 points $e_1,e_2,e_3,e_4,e_5$ in general
position in $\P{3}$ meet a plane $\pi$ not containing any of the points $e_i$
in a Desargues configuration. Conversely, it follows from the standard proof
of Desargues theorem via  $\P{3}$ (see \cite{kn:tod}) that every Desargues 
configuration is obtained
in this way. In the sequel, we choose the coordinates of $\p{3}$ 
in such a way that
\[e_1=(1:0:0:0),\dots,e_4=(0:0:0:1),e_5=(1:1:1:1).\]
If the plane $\pi$ is given by the equation  $\sum_{i=1}^4\alpha_ix_i=0,$
then the assumption that $e_i \not \in \pi$ for $i=1,\dots,5$ amounts
to 
\[ \alpha_i \neq 0 \mbox{ for } i=1,\dots,4 \mbox{ and } \sum_{i=1}^4 \alpha_i \neq 0. \]
We denote by $D_{\pi}$ the Desargues configuration determined by $\pi.$
It consists of the 10 points $p_{ij}:=\overline{e_ie_j} \cap \pi$  
and the 10 lines ${\ell}_{ijk}:=\overline{e_ie_je_k}\cap \pi$ for $i \leq i,j,k
\leq  5, i \neq j \neq k \neq i$ (see Figure 2). Note that the notation 
is meant to be symmetric, i.e. $p_{ij}=p_{ji}$ and 
$\ell_{ijk}=\ell_{jki},$ etc.
From the picture it is obvious that every point of a Desargues configuration
is the center of perspective of two triangles: the point $p_{ij}$ is the
centre of perspective of the triangles $p_{ik}p_{il}p_{im}$ and
$p_{jk}p_{jl}p_{jm}$ where $\{i,j,k,l,m \}=\{1,2,3,4,5\}.$
In particular every point of the configuration admits an axis of perspective: 
for $p_{ij}$ it is the line ${\ell}_{klm}.$

\begin{figure}
\hspace{2.5cm}
\caption{Desargues' Theorem in space}
\end{figure}

By definition two Desargues  configurations $D_1$ and $D_2$ are {\it isomorphic}
if there is an automorphism $\alpha \in PGL_2 (\complex)$ such that
$D_2=\alpha(D_1).$

\bigskip

\begin{lema} \label{carac}
For planes $\pi$ and $\pi' \in \p{3}$ not containing a point $e_i$
the following conditions are equivalent:

\noindent 1) $D_{\pi} $ is isomorphic to $D_{\pi'}.$ \newline
\noindent 2) There is an $A \in PGL_3(\complex)$ such that

{\rm a)} $A\pi=\pi',$

{\rm b)} A permutes the 5 points $e_1,\dots,e_5.$

\end{lema}

\bigskip

For the proof we need the following notion: 
a {\it complete quadrangle of a
Desargues configuration} $D$ is a set of 4 points and six
lines of $D$ which form the points and lines of a complete quadrangle.
From Figure 2 it is obvious that $D_{\pi}$ (and thus any $D$) admits 
exactly 5 complete quadrangles: any $e_i$  determines the complete quadrangle
consisting of the points $p_{ij},p_{ik},p_{il},p_{im}$ and the lines 
$\ell_{ijk},\ell_{ijl},\ell_{ijm},\ell_{ikl},\ell_{ikm},\ell_{ilm}.$

\bigskip

\textsc{ Proof of Lemma 1.1:} We have to show that (1) implies (2), the converse
implication being obvious. So let $\alpha:\pi \longrightarrow \pi' $
be a linear isomorphism with $\alpha(D_{\pi})=D_{\pi'}.$ We have to
show that $\alpha$ extends in a unique way to an $A \in PGL_3(\complex)$ permuting
the 5 points $e_1,\dots,e_5.$ Certainly $\alpha$ maps the 5 complete
quadrangles of $D_{\pi}$ onto the 5 complete quadrangles of $D_{\pi'}.$
As outlined above a complete quadrangle of $D_{\pi}$ is uniquely determined by
a point $e_i$ and similarly for $D_{\pi'}.$ Hence $\alpha$ induces a permutation
$\sigma$ of the 5 points $e_1,\dots,e_5.$ It is easy to see 
that there is a one-dimensional family 
$\{ A_t \in PGL_3(\complex)| t\in \P{1}\}$ satisfying $A_t(e_1)=e_{\sigma(1)}$
and $A_t(\pi)=\pi'$ (choose suitable coordinates for the source-$\P{3}$
and the image-$\P{3}$ of $A:\P{3} \longrightarrow \P{3}).$   
For every $t \in \P{1}$ the automorphism $A_t$ maps the line
$\overline{e_1p_{12}}$ onto the line 
$\overline{e_{\sigma(1)}p_{\sigma(1)\sigma(2)}}.$
Since $e_2 \in \overline{e_1p_{12}},$
$e_{\sigma(2)} \in \overline{e_{\sigma(1)}p_{\sigma(1)\sigma(2) }},$ and 
$\{ A_t(e_2) | t\in \P{1} \}=\overline{e_{\sigma(1)}p_{\sigma(1)\sigma(2)}},$
there is a unique
$t_0 \in \P{1}$ such that 
\[ A_{t_0} =e_{\sigma(2)}.\]
We claim that $A:=A_{t_0}$ satisfies apart from $(a)$ also condition $(b).$
But for $ 3 \leq k \leq 5$ we have $e_k=\overline{e_1p_{1k}}\cap 
\overline{e_2p_{2k}}.$
Hence 
\[A(e_k)=A(\overline{e_1p_{1k}})\cap A(\overline{e_2p_{2k}})=
\overline{e_{\sigma(1)}p_{\sigma(1)\sigma(k)}}\cap
\overline{e_{\sigma(2)}p_{\sigma_2\sigma(k)}}=e_{\sigma(k)} \]
This concludes the proof of the Lemma. \qed

\bigskip 

The 5 points $e_1,\dots,e_5$ form a projective basis of $\P{3}.$ Hence
for every permutation $\sigma$ of the points $e_1,\dots,e_5$ there is a unique
$A_{\sigma} \in PGL_3(\complex)$ inducing $\sigma.$ Let $\mathcal{ S}_5$ denote the
group of these $A_{\sigma}.$ The action of $\mathcal{ S}_5$ on $\P{3}$ induces
an action on the dual projective space $\PD{3}.$
Define 
\[ U:=\PD{3} \setminus \cup_{i=1}^5P_{e_i},\]
where $P_{e_i}$ denotes the plane in $\PD{3}$ parametrizing the planes $\pi$
with $e_i \in \pi.$ The action of $\mathcal{ S}_5$ on $\PD{3}$ restricts to an
action on $U.$
Since the quotient of a quasi-projective variety by a finite group is always
an algebraic variety, we obtain as an immediate consequence of 
Lemma~\ref{carac}

\bigskip

\begin{teo}
The algebraic variety $M_D:=U/ \mathcal{ S}_5$ is a moduli space
for Desargues configurations.
\end{teo}
 
\bigskip

It is clear how to define families of Desargues configurations.
Doing this it is easy to see that $M_D$ is a coarse moduli space in the sense
of geometric invariant theory. In particular, $M_D$ is uniquely determined 
as an algebraic variety. Since a Desargues configuration may admit a nontrivial 
group of automorphisms, $M_D$ is not a fine moduli space.

\bigskip

Next we work out the subspace of $M_D$ parametrizing special 
Desargues configurations. By definition a Desargues configuration $D$
 is called {\it special} if and only if $D$ contains a point lying on its axis. The axis of
a point $p_{ij}$ is the line $\ell_{klm}$ with $\{i,j,k,l,m\}=\{1,2,3,4,5\}.$
This implies that a configuration $D_{\pi}$ is special if and only if the
plane contains the point of intersection of a line $\overline{e_ie_j}$
with the plane $\overline{e_ke_le_m}.$
Hence a Desargues configuration $D_{\pi}$ is special if and only if
the plane $\pi$ contains a point $Q_{ij}=\overline{e_ie_j}\cap 
\overline{e_ke_le_m},\{i,j,k,l,m\}=\{1,\dots,5\}.$
Now let $P_{ij} \subset U = \PD{3} \setminus \cup_{i=1}^5P_{e_i},$
denote the hyperplane parametrizing the planes $\pi$ in $\P{3}$
not containing the points $e_i,$ but containing the point $Q_{ij}.$
If $q:U \longrightarrow M_{D}$ denotes the natural projection 
map, then we conclude, since the group $\mathcal{ S}_5$ obviously acts transitively
on the set of points $\{Q_{ij}\}:$

\bigskip

\begin{prop}
The special Desargues configurations are parametrized  by the irreducible
divisor $q(P_{12})$ in $M_D.$
\end{prop}

\bigskip

In other words: if $D$ is a special Desargues configuration then there is
a plane $\pi \subset \P{3}$ containing the point $P_{12}=(1:1:0:0)$ and not
containing any $e_i$ such that $D_{\pi} \cong  D.$

\bigskip

\bigskip
\begin{figure}
\hspace{4cm}
\caption{A degenerate Desargues configuration of the 1st kind}
\end{figure}

\bigskip

Finally we introduce degenerate Desargues configurations.
Desargues' Theorem fails if one of the two triangles is replaced by three
lines passing through the center of perspective. However, the following
classical theorem (see \cite{kn:SK}, p.99, Exercise 14) may be considered as
a limiting case of Desargues' Theorem: If $ABC$ is a proper triangle and $a,b,c$
are 3 lines passing through a point $O$ (not lying on the sides of the
triangle) then the points $a \cap \overline{BC},b \cap \overline{AC}, c \cap
\overline{AB}$ are collinear if and only if there is a an involution $\imath$
of the $\P{1}$ of lines centered in $O$ such that   $\imath(a)=\overline{OA},
\imath(b)=\overline{OB},\imath(C)=\overline{OC}.$

It is now easy to see that for any configuration of 7 points and 10 lines
in $\P{2}$ satisfying the conditions of the theorem there is a plane
$\pi \subset \P{3}$ containing exactly 1 of the 5 points $e_i$ such that the
configuration is isomorphic to the configuration $D$ cut out on $\pi$ by the
lines $\overline{e_ie_j}$ and the planes $\overline{e_ie_je_k}.$
A picture of a degenerate Desargues configuration $D_{\pi}$
is shown in Figure 3 where $e_5 \in \pi$ and hence $p_{15}=p_{25}=p_{35}= 
p_{45}=e_5$ with the notation above.


\bigskip

Note that of the 5 complete quadrangles of the Desargues configuration
4 still survive in the degenerate case namely (in the case of Figure 3)
$\{p_{12},p_{13},p_{14},p_{15}\}$,$\{p_{12},p_{23},p_{24},p_{25}\}$,
$\{p_{13}.p_{23},p_{34},p_{35}\}$,$\{ p_{14},p_{24},p_{34},p_{45} \}.$
This implies that the proof of Lemma 1.1 also works in this case. In other
words: two degenerate Desargues configurations $D_{\pi}$ and $D_{\pi'}$
are isomorphic if and only if there is an automorphism $A \in PGL_3(\complex)$
with $A\pi=\pi'$ permuting the 5 points $e_1,\dots,e_5$ in $\P{3}.$
Similar remarks can be made if the configuration is even more degenerate, that
is if 2 triangles of a Desargues configuration collapse (equivalently
if the plane $\pi$ passes through 2 of the points $e_i$ (see Figure 4 where
the triangle $p_{14}p_{24}p_{34}$  of Figure 3 collapsed) or if 2 triangles
collapse and the 2 points come together (equivalently if the plane $\pi$ passes
through 3 of the points $e_i$). In this case the configuration is a complete
quadrangle (see Figure 5, where $e_3=p_{13}=p_{23})$. We omit the details
since they are easy to work out (see also Section 6).
\bigskip

\bigskip

These  remarks induce the following definitions: as above fix 5 points
$e_1,\dots,e_5 \in \P{3}$ in general position and consider \textbf{ any }
 plane 
$\pi \subset \P{3}.$ The 10 lines $\overline{e_ie_j}$ and 10 planes 
$\overline{e_ie_je_k}$ in $\P{3}$ cut out on $\pi$ a configuration $D_{\pi}$
of points and lines. We call $D_{\pi}$ a {\it generalized Desargues
configuration.} If $\pi$ contains a point $e_i,$ $D_{\pi}$ is called 
a {\it degenerate
Desargues configuration.} It is called {\it of the $i$-th kind} if $\pi$ contains exactly 
$i$ of the points $e_i$ for $i=1,2,3.$

\bigskip

It is easy to see that Lemma 1.1 remains valid for generalized Desargues
configurations. In fact the same proof works also in the degenerate case.
One has only to remark that a degenerate Desargues configuration $D_{\pi}$
of the first kind (respectively $2^{nd}$ kind, respectively $3^{rd}$ kind)
admits 4 (respectively 3, respectively 1) complete quadrangles.
In the same way as we deduced Theorem 1.2 from Lemma 1.1 we obtain from this

\bigskip

\begin{teo}
The variety $\overline{M}_D:=\PD{3}/ \mathcal{ S}_5$ is a moduli space for generalized
 Desargues configurations.
 
\end{teo}

\begin{rmk}
{\rm In \cite{kn:Ma} which is the standard reference for Desargues
configurations, K. Mayer constructs the moduli space $M_{D}$ in a different
way. He chooses the 4 points of a complete quadrangle as a projective
basis of $\P{2},$ say $A,A_1,B_1,C_1$ in Figure 1.
Consider the points $A_1':=\overline{AA_1}\cap \overline{B_1C_1},
B_1':=\overline{AB_1}\cap\overline{A_1C_1}$ and $ C_1'=\overline{AC_1}\cap
\overline{A_1B_1}.$
The three crossratios $(A,A_1,A_1',A_2), (A,B_1,B_1',B_2)$ and
$(A,C_1,C_1',C_2)$ determine the Desargues configuration and an open set
$V$ of $\P{1}^3$  represents Desargues configurations. The choice of
a complete quadrangle induces an action of $\mathcal{ S}_5$  on $V$ which
Mayer worked out explicitly. Moreover he showed that 
$V / \mathcal{ S}_5$
is the moduli space of Desargues configurations. For our purposes, the
construction of Theorem 1.2 is more appropriate.} 
\end{rmk}  

\bigskip

\begin{figure}
\hspace{4cm}
\caption{A degenerate Desargues configuration of the 2nd kind}
\end{figure}

\section{The map $\Phi:M_D \longrightarrow \sex$}

Recall that  a  binary form  $f(x,y)$ of degree 6 is
called {\it stable} if f admits no root of multiplicity $\geq 3$ and that
the space $\sex$ of stable binary sextics exists (see \cite{kn:MF}).
In this section, we present the construction of Stephanos (see \cite{kn:ST}
-in a slightly different set up) associating to every Desargues configuration
$D_{\pi}$ a binary sextic $J_{\pi}$ in a canonical way. In Section~\ref{stable}
we will see that $J_{\pi}$ is stable so that we get a holomorphic map 
$\Phi: M_D \longrightarrow \sex.$

\bigskip

As in the last section we fix the coordinate system of $\P{3}$ in such a way
that $e_1=(1:0:0:0),\dots,e_4=(0:0:0:1)$ and $e_5=(1:1:1:1).$
Let $\pi$ be a {\it nondegenerate plane} in $\P{3},$ meaning that the configuration
$D_{\pi}$ is nondegenerate, with equation $\sum_{i=1}^4\alpha_iz_i=0.$
For the coefficients $\alpha_i$ of $\pi$ this just  means $\alpha_i \neq 0$
for $i=1,\dots,4$ and $\sum_{i=1}^4 \alpha_i \neq 0.$   
In 1847, von Staudt proved in his fundamental book \cite{kn:VS} that there is a 
unique (smooth) conic $s_{\pi}$ on $\pi$ such that the polar line of every point
of the configuration $D_{\pi}$ is its axis. We call $s_{\pi}$ the {\it von Staudt
conic } of $D_{\pi}.$

\bigskip

\begin{lema} \label{VS}
The von Staudt conic $s_{\pi}$ of $D_{\pi}$ is given by the equation 
$\sum_{i=1}^4\alpha_iz_i^2=0$ in $\pi:\sum_{i=1}^4\alpha_iz_i=0$
\end{lema}

\bigskip

\textsc{ Proof:} According to a remark of Reye, which is easy to check (see
\cite{kn:Re},p. 135), there is a unique quadric $S_5$ in $\P{3}$ such that 
the tetrahedron $e_1,e_2,e_3,e_4$ is a polar tetrahedron, i.e. the polar plane
 of $e_i$ with respect to $S_5$ is the opposite plane of the tetrahedron, and
such that the polar plane of $e_5$ is the plane $\pi.$
If, as usual, we denote by $S_5$ also the matrix of the quadric, the
conditions mean $e_iS_5e_j=0$ for $1 \leq i,j \leq 4,$ $i\neq j$ and
$e_5S_5=(\alpha_1,\alpha_2,\alpha_3,\alpha_4)^t.$ But this implies
$S_5=diag(\alpha_1,\alpha_2,\alpha_3,\alpha_4).$
Comparing the definitions, the von Staudt conic $s_{\pi}$ is 
just the restriction of $S_5$ to $\pi,$ which gives the assertion. \qed
\bigskip

Note that the proof of Lemma~\ref{VS} 
yields actually a proof of von Staudt's Theorem. Obviously $S_5$ depends
on the choice of $e_5$ as the cone vertex. One could in an analogous way
define $S_i \; (1\leq i \leq 4)$ and use it to prove Lemma~\ref{VS}.

\bigskip

As remarked in the last section, each point $e_i$ corresponds
to a complete quadrangle of $D_{\pi},$ namely
$p_{1i},\dots,\check{p_{ii}},\dots, p_{5i}.$
Since these 4 points are in general position, there is a unique pencil of
conics $\lambda q_1+\mu q_2$ passing through the 4 points. The pencil
$\lambda q_1+\mu q_2$ restricts to a pencil of effective divisors of degree
4 on the von Staudt conic $s_{\pi}.$ The Jacobian of this pencil
is an effective divisor $j_{\pi}$ of degree 6
on $s_{\pi}$. It is defined as follows: choose
an isomorphism $\varphi:s_{\pi}\longrightarrow \P{1}$ and coordinates
$(x_1,x_2)$ of $\P{1}.$
There are binary quartics $f_1(x_1,x_2)$ and $f_2(x_1,x_2)$
with zero divisor $\varphi(q_1\cap s_{\pi})$ and 
$\varphi(q_2\cap s_{\pi}).$
The Jacobian of $f_1$ and $f_2$

\[ J_{\pi}(x_1,x_2)=det(\frac{\partial f_i}{\partial x_j}) \]
is a binary sextic with zero divisor $\varphi(j_{\pi}).$
Note the binary sextic is determined by the pencil $\lambda q_1+\mu q_2$
up to a nonzero constant and an automorphism of $\P{1}.$ We call $J_{\pi}$
the {\it Jacobian sextic} of the pencil, for the Jacobian divisor $j_{\pi}$
on $s_{\pi}.$ 

\bigskip

\begin{teo} \label{jac}
The Jacobian divisor $j_{\pi}$ on the von Staudt conic $s_{\pi}$
is uniquely determined by the Desargues configuration $D_{\pi}.$
It is cut out on $s_{\pi}$ by the cubic $\sum_{i=1}^4\alpha_iz_i^3=0.$
\end{teo}

\bigskip

\textsc{ Proof:} Choose first $i=5.$  The pencil of conics passing through  the 4
points $p_{15},p_{25},p_{35},p_{45}$ is cut out on $\pi$ by the pencil
of cones in $\P{3}$ passing through the 4 lines $\overline{e_1e_5},
\overline{e_2e_5},\overline{e_3e_5},\overline{e_4e_5}.$
The conditions for a quadric $Q \subset \p{3}$ to be in this pencil are:
$Q$ is singular in $e_5,$ i.e. $Qe_5=0$ and $e_i \in Q$ for $i=1,2,3,4,$
i.e. $e_iQe_i=0.$
This yields
\[ Q= \left( \begin{array}{cccc}
         0&a&b&-(a+b)\\
         a&0&-(a+b)&b\\
         b&-(a+b)&0&a\\
         -(a+b)&b&a&0

\end{array}\right) \]
Hence the pencil is generated by the cones 
\[ Q_1= \left( \begin{array}{cccc}
         0&1&0&-1\\
         1&0&-1&0\\
         0&-1&0&1\\
         -1&0&1&0

\end{array}\right)  \mbox{ and }
 Q_2= \left( \begin{array}{cccc}
         0&0&1&-1\\
         0&0&-1&1\\
         1&-1&0&0\\
         -1&1&0&0

\end{array}\right). \]

It is easy to see that the Jacobi divisor of the pencil $\lambda q_1+\mu q_2$
on the curve $s_{\pi}=\{
(z_1:\dots:z_4) \in \P{3} |
\sum_{i=1}^4\alpha_iz_i=\sum_{i=1}^4\alpha_iz_i^2=0 \}$
is the intersection of the curve $s_{\pi}$ with the hypersurface of points in 
$\P{3}$ where the vectors of partial derivatives $d(\sum\alpha_iz_i)$, $d(\sum
\alpha_iz_i^2)$, $d(Q_1)$, $d(Q_2)$ are linearly dependent.
Here we consider $Q_i$ as the quadratic form determined by the matrix $Q_i.$
This hypersurface is given by the equation:

\[ J:= det \left( \begin{array}{cccc}
         \alpha_1&\alpha_1z_1&z_2-z_4&z_3-z_4\\
\alpha_2&\alpha_2z_2&z_1-z_3&z_4-z_3\\
 \alpha_3&\alpha_3z_3&z_4-z_2&z_1-z_2\\        
\alpha_4&\alpha_4z_4&z_3-z_1&z_2-z_1
\end{array}\right) =0. \]
We claim that $J$ is congruent to $\sum_{i=1}^4\alpha_iz_i^3$ modulo
the ideal $(\sum_{i}\alpha_iz_i,\sum_{i}\alpha_iz_i^2).$
But
\begin{eqnarray*}
 J&=& det \left( \begin{array}{cccc}
         \alpha_1&\alpha_1z_1&z_2-z_4&z_3-z_4\\
\alpha_2&\alpha_2z_2&z_1-z_3&z_4-z_3\\
 \alpha_3&\alpha_3z_3&z_4-z_2&z_1-z_2\\        
\sum \alpha_i&\sum \alpha_iz_i&0&0
\end{array}\right)\\
& \equiv& -\sum_i\alpha_i det \left( \begin{array}{ccc}
         \alpha_1z_1&z_2-z_4&z_3-z_4\\
         \alpha_2z_2&z_1-z_3&z_4-z_3\\
         \alpha_3z_3&z_4-z_2&z_1-z_2
\end{array}\right) \pmod{\sum_{i} \alpha_iz_i} \\
&=&-\sum_{i}\alpha_i [2(\alpha_1z_1^3+\alpha_2z_2^3+\alpha_3z_3^3)
-(\alpha_1z_1^2-\alpha_2z_2^2-\alpha_3z_3^2)z_4
-z_4(\alpha_1z_1^2+\alpha_2z_2^2+\alpha_3z_3^2)\\
&&+(z_1+z_2+z_3)(z_4(\alpha_1z_1+\alpha_2z_2+\alpha_3z_3)
-(\alpha_1z_1^2+\alpha_2z_2^2+\alpha_3z_3^2))] \\
& \equiv& -2 \sum_{i=1}^4 \alpha_i \sum_{j=1}^4 \alpha_jz_j^3 \;\;
\pmod{(\sum_{i}\alpha_iz_i, \sum_{i}\alpha_iz_i^2)}
\end{eqnarray*}

It remains to show that if we start with another point $e_i$
we obtain the same Jacobi divisor $j_{\pi}$ on $s_{\pi}.$ Without loss of generality
we may start with $e_1,$  since choosing another $e_i \;(2 \leq i \leq 4)$
amounts only to a permutation of the coordinates.

The pencil of conics in $\pi$ passing through the points $p_{12},
p_{13},p_{14},p_{15}$ is cut out on $\pi$ by the pencil of cones in $\P{3}$ 
passing through the lines $\overline{e_1e_2},\overline{e_1e_3},
\overline{e_1e_4},\overline{e_1e_5}.$
This pencil is generated by the cones:
\[ Q_1'= \left( \begin{array}{cccc}
         0&0&0&0\\
         0&0&1&0\\
         0&1&0&-1\\
         0&0&-1&0

\end{array}\right)  \mbox{ and }
 Q_2'= \left( \begin{array}{cccc}
         0&0&0&0\\
         0&0&0&1\\
         0&0&0&-1\\
         0&1&-1&0

\end{array}\right). \]
A very similar computation as the one  above shows that modulo the ideal 
$(\sum \alpha_i z_i,\sum \alpha_iz_i^2)$ the hypersurface of points in $\P{3}$
where the vectors of partial derivatives $d(\sum \alpha_iz_i)$, $d(\sum
\alpha_iz_i^2)$,
 $d(Q_1')$, $d(Q_2')$ are linearly dependent is given by 
$\sum _{i=1}^4 \alpha_i z_i^3,$ which completes the proof of the theorem.
\qed

\bigskip

If $\pi'$ is another nondegenerate plane in $\P{3}$ such that the Desargues
configurations $D_{\pi}$ and $D_{\pi'}$ are isomorphic, it follows immediately
from Lemma~\ref{carac} and Theorem~\ref{jac} that the isomorphism 
$\alpha: \pi \longrightarrow \pi'$ with $\alpha(D_{\pi})=D_{\pi'}$ maps the
Jacobi divisor $j_{\pi}$ on $s_{\pi}$ onto the Jacobi divisor $j_{\pi'}$ on
$s_{\pi'}.$ 

In the sequel we will always interpret $j_{\pi}$ as the binary sextic $J_{\pi}$
(see above and recall that $J_{\pi}$ is determined up to a nonzero constant
and up to an automorphism of $\P{1}$).
In Theorem~\ref{sta}, we will show that $J_{\pi}$ is always stable.
Hence we obtain a canonical map 
\[ \Phi: M_D \longrightarrow \sex \]
of the moduli space of Desargues configurations into the moduli space
of stable binary sextics. It is clear that $\Phi$ is holomorphic, since
Theorem~\ref{jac} implies that it is given by polynomials.

\bigskip

\begin{rmk}
{\rm The proof of Theorem~\ref{jac} suggests  another definition of the Jacobi
divisor $j_{\pi}:$ the von Staudt conic $s_{\pi}$ is not contained in the
pencil of conics $\lambda q_1+\mu q_2$ passing through the points
$p_{15},p_{25},p_{35},p_{45}.$ Hence $\lambda q_1+\mu q_2+\nu s_{\pi}$
is a net of conics in the plane $\pi.$ Its discriminant locus
 given by
\[ det(q_1,q_2,s_{\pi})=0\]
is a cubic which intersects the conic $s_{\pi}$ in the Jacobian divisor
$j_{\pi}.$ For the proof note only that the net is general and thus its
discriminant is equivalent to its Jacobian locus given by $det(dq_1,dq_2,ds_{\pi})=0.$}  
\end{rmk}   

\begin{rmk} \label{config}
{\rm Given a nondegenerate plane $\pi,$ we associated to every point $e_i,
(1 \leq i \leq 5)$ a pencil of quartic divisors $\lambda f_1^i + \mu f_2^i$ on
the von Staudt conic $s_{\pi}.$ It is easy to check that in general these
5 pencils are different from each other. }

\end{rmk}

\section{Injectivity of $\Phi:M_D \longrightarrow \sex$  }

The main object of this section is to prove the following theorem.

\bigskip

\begin{teo} \label{iso}
The map $\Phi: M_D \longrightarrow \sex$ is an injective birrational morphism.
\end{teo}

\bigskip

For the proof we need some preliminaries. Let $A$ and $B$ denote two quadrics
in $\P{n}$ (for us $n=2$ or $3)$ 
given by the equations $x^tAx=0$ and $x^tBx=0.$
Recall that $A$ is called {\it apolar} to $B$ if
\[ tr(A \cdot Adj(B))=0.\]
Here $adj(B)$ denotes the adjoint matrix of $B,$ i.e. the matrix of the dual
quadric $\hat{B}$ of $B.$ 
Note that this definition is not symmetric (some authors say that $A$ is
apolar to $\hat{B}$).
Geometrically this means the following (see \cite{kn:SK}):
a conic $A$ is apolar to $B$ in $\P{2}$ if and only if there is a triangle
inscribed in $A$ and self-polar with respect to $B.$ A quadric $A$ 
is apolar to $B$
in $\P{3}$ if and only if there is a tetrahedron inscribed in $A$
and self-polar with respect to $B.$

\begin{prop} \label{apol}
Let $D_{\pi}$ be a Desargues configuration. Every conic in $\pi$  passing
through the 4 points of a complete quadrangle of $D_{\pi}$ is apolar to the von
Staudt conic $s_{\pi}.$ 
\end{prop}

\bigskip

\textsc{Proof:}
Since apolarity is a linear condition and by the special choice of the
coordinates it suffices to show that the quadric cones $Q_1,Q_2,Q_1'
Q_2'$ (for the notation see proof of Theorem~\ref{jac}) 
are apolar to the the quadric $diag(\alpha_1,\alpha_2,\alpha_3,\alpha_4)$
in $\P{3},$ whose restriction to $\pi$ is $s_{\pi},$ which is an immediate
computation.  \qed

\bigskip

In order to show that a conic $q$ in the plane $\pi$ passing through the
points of a complete quadrangle of $D_{\pi}$ is the unique conic apolar to
$s_{\pi}$ and passing through $q \cap s_{\pi}$(with multiplicities) we change
the coordinates. We choose the coordinates $(x_0:x_1:x_2)$ of $\pi=\P{2}$
in such a way that the von Staudt conic $s_{\pi}$ is given by the equation
\begin{eqnarray} \label{con}
 x_1^2-4x_0x_2=0
\end{eqnarray}

So the matrix of the dual conic is:  
\[ \widehat{s_{\pi}}= \left( \begin{array}{ccc}
         0&0&2\\
         0&-4&0\\
         2&0&0

\end{array}\right)  \]
Let $(t_0:t_1)$ denote homogeneous coordinates of $\P{1}.$ Then
\[ (t_0:t_1) \longrightarrow (x_0:x_1:x_2)=(t_0^2: 2t_0t_1: t_1^2) \]
is a parametrization $\P{1}\stackrel{\sim}\longrightarrow s_{\pi}$ of
$s_{\pi}.$ Using this, any effective divisor of degree $4$ on $s_{\pi}$
can be interpreted as a binary quartic $f(t_0,t_1).$ Hence
\begin{eqnarray} \label{qua}
f(t_0,t_1)&=&a_0t_0^4+4a_1t_0^3t_1+6a_2t_0^2t_1^2+4a_3t_0t_1^3+a_4t_1^4 
\end{eqnarray}
with $(a_0:a_1:a_2:a_3:a_4) \in \P{5}$ representing an arbitrary divisor
of degree 4 on $s_{\pi}.$ Given an effective divisor $\delta,$ of degree 4 on
$s_{\pi},$ i.e. a binary quartic  
 $f,$ consider the conic
\[ q=a_0x_0^2+a_2x_1^2+a_4x_2^2+2a_1x_0x_1+2a_2x_0x_2+2a_3x_1x_2. \]

\begin{lema} \label{ass}
The conic $q$ is the unique conic passing through the divisor $\delta$
and apolar to the von Staudt conic $s_{\pi}.$
\end{lema}
  
\textsc{Proof:}
Note first that $q$ passes through $\delta$ 
by the choice of the parametrization. It is apolar to $s_{\pi},$ since

\[ tr \left( \left( \begin{array}{ccc}
         a_0&a_1&a_2\\
         a_1&a_2&a_3\\
         a_2&a_3&a_4
\end{array}\right) 
 \left( \begin{array}{ccc}
         0&0&2\\
         0&-4&0\\
         2&0&0
\end{array}\right) \right)=2a_2-4a_2+2a_2=0
. \]
On the other hand, apolarity is one linear condition. Hence the
space of conics apolar to $s_{\pi}$ is isomorphic to $\P{4}.$
Now it is easy to see (\cite{kn:CX}, p.17) that any divisor of degree 4 on $s_{\pi}$ imposes
independent conditions on conics . In other words, the space of conics
passing through $\delta$ is isomorphic to $\P{1}.$ Hence there is exactly one
conic passing through $\delta$ and apolar to $s_{\pi}$ if we only show that
not every conic passing through $\delta$ is  apolar to $s_{\pi}.$

For this we assume $\delta= p_1+p_2+p_3+p_4$ with $p_i \neq p_j,$
for $i \neq j,$ the degenerate cases being even easier to check.
We may choose the coordinates in such a way that in $\P{1}$ we have 
$p_1=(0:1),$ $p_2=(1:0),$ $p_3=(1:1),$ $p_4=(1:t)$ with $t \neq 0,1.$
So in $\P{2}:$
$p_1=(0:0:1),p_2=(1:0:0),p_3=(1:2:1),p_4=(1:2t:t^2)$
then $2t(t-1)xy+(1-t^2)y^2+(2t-2)yz=0$ is a conic passing through
$p_i,$ for $i=1,\dots,4$ and not apolar to $s_{\pi},$ if $t\neq -1.$
\newline\noindent
If $t=-1$ then $2xy+4xz-y^2-2yz$ satisfies these conditions.
\qed

\bigskip

We call $q$ the conic {\it associated to the binary quartic}~(\ref{qua}).
If we associate to every divisor of a pencil $\lambda \delta_1 + \mu \delta_2$
of quartic divisors the unique conic of Lemma~\ref{ass},
we obtain a pencil of conics. Hence we obtain as  an immediate consequence
of Proposition~\ref{apol} and Lemma~\ref{ass}:

\begin{cor}
Let $D_{\pi}$ be a Desargues configuration.
The pencil of conics $\lambda q_1 + \mu q_2$ passing through the 4 points of a
 complete quadrangle of $D_{\pi}$ is the unique pencil of conics cutting out
the pencil of quartic divisors $(\lambda q_1+\mu q_2) \cap s_{\pi}$ which is
apolar to the von Staudt conic $s_{\pi}.$
\end{cor}

\bigskip

Given a smooth conic $q_0,$ we choose the coordinates in such a way that
$q_0$ is given by equation (\ref{con}).
As outlined above, any pencil of binary quartics can be interpreted
(up to isomorphism) as a pencil of quartic divisors on $q_0.$
There is a unique pencil of conics associated to it according 
to Lemma~\ref{ass}. We call a pencil of binary quartics
(respectively the corresponding pencil of quartic divisors) {\it admissible}
if its associated pencil of conics is general, i.e. its base locus
consists of 4 different points. The following proposition is the
first step in the proof of Theorem~\ref{iso}.

\bigskip

\begin{prop}\label{admiss}
Consider an admissible pencil of quartic divisors on a smooth conic $q_0.$
There is a unique Desargues configuration $D$ such that $q_0$ is the von
Staudt conic of $D$ and the pencil of quartics is cut out on $q_0$ 
by the pencil of conics  passing through the points of a complete
quadrangle of $D.$ 
\end{prop} 
\textsc{Proof:}
Let $p_1,p_2,p_3,p_4$ denote the 4 base points of the pencil of conics 
associated to the given pencil of quartic divisors.
Consider the complete quadrangle consisting of the 4 points $p_1,p_2,p_3,p_4$
and the 6 lines $\ell_{ij}=\overline{p_ip_j}$ for $1 \leq i < j \leq 4.$
Let $p_{ij}$ be the poles of the lines $\ell_{ij}$ and $\ell_i$ the polars 
of the points $p_i$ with respect to the conic $q_0.$
Then the 10 points  $P_{ij},p_i$ and the 10 lines $\ell_{ij},\ell_i$ form 
a Desargues configuration according to a theorem proved by von Staudt  
in 1831 (see \cite{kn:BLA}, p.62). Using Lemma~\ref{ass} we have the assertion.
\qed

\bigskip

Before we go on, let us note the following characterization of admissible
binary quartics, which we need later.

\bigskip
\begin{prop} \label{admiss1}
The pencil of binary quartics $\lambda f_1 + \mu f_2$ with 
$f_i=\sum_{j=1}^4 {4 \choose i} a_j^it_0^{4-j}t_1^j$ is admissible if and only if the 
discriminant of the binary cubic $det(t_0 q_1+t_1q_2)$ is nonzero, where
$ q_i= \left( \begin{array}{ccc}
         a_0^i&a_1^i&a_2^i\\
         a_1^i&a_2^i&a_3^i\\
         a_2^i&a_3^i&a_4^i
\end{array}\right)$
is the matrix of the conic associated to $f_i$
for $i=1,2.$
  
\end{prop}
\bigskip

\textsc{Proof:}
This is a consequence of the fact (see \cite{kn:tod}, Section 6.3,)
that the base locus of a pencil of conics consists of 4 different points
if and only if its discriminant does not vanish
together with Lemma~\ref{ass}.\qed

\bigskip
In order to complete the proof of Theorem~\ref{iso}, we need to compute 
the degree of the Jacobian map. Note that the  
binary quartic~(\ref{qua}) determines a point in $\P{4}$ namely 
$(a_0:\dots:a_4).$
If we consider $\P{4}$ as the space of quartics, the space of pencils of
binary quartics is $Gr(1,4)$ the Grassmannian of lines in $\P{4}.$
Considering in the same way the space of binary sextics
$\sum_{i=0}^6a_it_0^{6-i}t_1^i$ as $\P{6},$  the map associating to every
pencil of quartics its Jacobian defines a morphism 
\[ Jac: Gr(1,4) \longrightarrow \P{6}. \]
Note that $Gr(1,4)$ is also of dimension 6. We need the following lemma:
\begin{lema} \label{deg5}
The Jacobian map $Jac : Gr(1,4) \longrightarrow \P{6}$ is a ramified
covering of degree 5.
\end{lema}

\textsc{Proof:}
Let $\lambda f_1+ \mu f_2$ be a pencil of quartics with 
$f_1=\sum_{i=1}^4 {4 \choose i} a_it_0^{4-i}t_1^i$ and
$f_2=\sum_{i=1}^4 {4 \choose i} b_it_0^{4-i}t_1^i.$ Its Jacobian is given by:

\[Jac(f_1,f_2)=det(\frac{\partial f_i}{\partial t_j})=\]
\[=16[
p_{01}t_0^6+3p_{02}t_0^5t_1+3(2p_{12}+p_{03})t_0^4t_1^2+
(8p_{13}+p_{04})t_0^3t_1^3+ 3(2p_{23}+p_{14})t_0^2t_1^4+
3p_{24}t_0t_1^5+p_{34}t_1^6] ,\]
where $p_{ij}:=a_ib_j-a_jb_i$ denotes the Pl\"ucker coordinates
of the pencil.
In particular up to a constant the Jacobian  does not depend on
the choice of $f_1$ and $f_2$ and defines a point $Jac(\lambda f_1+\mu f_2)$
in $\P{6}.$
The explicit form of $Jac(f_1,f_2)$ implies that the map $Jac$ 
factorizes via the Pl\"ucker embedding $p$
\bigskip
\begin{center}

\xymatrix{ &&&&&&&\P{9}\ar@{.>}[d]_q\\
           &&&&&Gr(1,4)\ar@{^{(}->}[urr]^p \ar[rr]^{Jac}&&\P{6}
}

\end{center}

\bigskip
with a linear projection map $q.$
In fact $Jac$ is given by
\[ \lambda f_1+\mu f_2 \mapsto (p_{01}:3p_{02}:3(2p_{12}+p_{03}):
8p_{13}+p_{04}:3(2p_{23}+p_{14}):3p_{24}:p_{34}) \in \P{6} .\]
It is the linear projection of the $\P{9}$ with coordinates 
$(p_{01}:\dots:p_{34})$ with center the plane $P$ with equations
$p_{01}=p_{02}=2p_{12}+p_{03}=8p_{13}+p_{04}=2p_{23}+p_{14}=p_{24}=p_{34}=0.$
On the other hand, the Pl\"ucker variety $p(Gr(1,4))$ in $\P{9}$ is given by the
Pl\"ucker relations among which, there are
\begin{eqnarray*}
p_{12}p_{34}-p_{13}p_{24}+p_{14}p_{23}&=&0\\    
p_{01}p_{34}-p_{03}p_{14}+p_{04}p_{13}&=&0\\     
p_{01}p_{23}-p_{02}p_{13}+p_{03}p_{12}&=&0        
\end{eqnarray*}

Now if 
$(p_{01}:\dots:p_{34}) \in P\cap p(Gr(1,4)),$ this implies 
$p_{14}p_{23}=p_{04}p_{13}=p_{03}p_{12}=0,$ so all $p_{ij}=0.$
Hence the center of projection $P$ does no intersect the
Pl\"ucker variety.
This implies that the degree of $Jac$ equals the degree of the Pl\"ucker variety
 $p(Gr(1,4))$ in $\P{9}.$
It is well known that this degree is 5 
(see \cite{kn:JH},p. 247). Since $Jac$ cannot contract a positive
dimensional subvariety of $Gr(1,4)$ this completes the proof of the lemma.
\qed

\bigskip

Using this we can finally prove Theorem~\ref{iso}.

\bigskip

\textsc{Proof of Theorem~\ref{iso}:}
Since admissibility of a pencil of quartics is an open condition and the
Jacobian map is finite, we conclude from Proposition~\ref{admiss}
that $\Phi$ is dominant  and there is an open dense set 
$U\subset \sex$ such that for every sextic $j \in U$
the preimage $Jac^{-1}(j)$ consists of exactly 5 pencils of quartics.
Shrinking $U$ if necessary we may assume according to 
Remark~\ref{config}
that every Desargues configuration $D \in \Phi^{-1}(U)$ admits 5
different pencils of quartic divisors on the von Staudt conic. Hence 
for every $j \in U$ the 5 pencils of quartics of any preimage
$D \in \Phi^{-1}(j)$ must coincide. Thus by Proposition~\ref{admiss}
 $\Phi^{-1}(j)$ consists only of one Desargues configuration, i.e.
$\Phi: \Phi^{-1}(U) \longrightarrow U$ is bijective. This implies 
that $\Phi$ is of degree 1. It remains to show the injectivity.

Since $\Phi$ is a finite birrational morphism of normal varieties, the
injectivity can only fail at points $p$ of $M_D$ for which $\Phi(p)$ is a
singular point of $\sex.$
But $\sex$ admits only 3 singular points (see \cite{kn:pqh}, Remark 5.4,
p.5551 ),
namely a point of mutiplicity 5 represented by the binary sextic $j_5=t_0^6
-t_0t_1^5,$ a point of multiplicity 2 represented by 
$j_2=t_0^6-t_0^2t_1^4,$
 a point of multiplicity 3 represented by $j_3=t_0^5t_1-t_0^2t_1^4.$
The moduli space $M_D$ also admits a singular  point  of multiplicity 5:
the Desargues configuration $D_{{\pi}_5}$ with $\pi_5$ given by 
$z_1-\mu z_2+\mu^2 z_3- \mu ^3z_4=0$ where $\mu$ is a root of the 
equation
$x^4+x^3-x^2+x-1=0,$ admits an automorphism of order 5 namely:  
\[ \left( \begin{array}{cccc}
         0&0&0&1\\
         -1&0&0&1\\
         0&-1&0&1\\
         0&0&-1&1 
\end{array}\right).\] It permutes the points $e_i$ as follows:
$e_1 \mapsto e_2 \mapsto e_3 \mapsto e_4 \mapsto e_5 \mapsto e_1.$
It is easy to check that $D_{{\pi}_5}$ is of multiplicity 5 in $M_d$ 
and that 
$\Phi^{-1}(j_5)=D_{{\pi}_5}.$

As in Example~\ref{j2} below one checks that all pencils of quartics
with Jacobian $j_2$ are not admissible. Hence $j_2$ is not in the image of $\Phi.$
Finally as in Example~\ref{j3} below one checks that only one pencil of
quartics with Jacobian $j_3$ is admissible namely $\lambda f_3+\mu g_3$ with
$f_3=2t_0^4+6t_0^2t_1^2+4t_0t_1^3$ and $g_3=\frac{1}{16} t_0^2t_1^2.$ Hence
$\Phi$ is also injective at the corresponding point at $M_D.$ 
This completes the
 proof of the theorem. \qed

\section{The image of the map $\Phi$} \label{stable}

Using the set up of the last section we are now in a position to prove the
following theorem, which was anounced (but not applied) already in
Section 3.

\bigskip

\begin{teo} \label{sta}
For any Desargues configuration $D$ the associated binary sextic 
$J=\Phi(D)$
is stable.
\end{teo}  

\bigskip

\textsc{Proof:}
Let $J(t_0,t_1)$ be a nonstable binary sextic, i.e admitting a root of multiplicity
$\geq 3.$ We choose coordinates in such a way that this root is
$(0:1).$ Hence $J$ is of the form
\[J(t_0,t_1)=A_0t_0^6+A_1t_0^5t_1+A_2t_0^4t_1^2+A_3t_0^3t_1^3 .\]
Suppose $\lambda f+\mu g$ with 
$f=\sum_{i=1}^4 {4 \choose i} a_it_0^{4-i}t_1^i$ and
$g=\sum_{i=1}^4 {4 \choose i} b_it_0^{4-i}t_1^i$
is a pencil of binary quartics whose Jacobian is $J$ (up to a nonzero
constant). According to the results of Section 4, it suffices to show
that $\lambda f + \mu g$ is not admissible which according to 
Proposition~\ref{admiss1} means that $det(t_0q_1+t_1q_2)$ admits
a double root where $q_1$ and $q_2$ denote the conics associated 
to $f$ and $g.$  
By the equation for $J$ given in the proof of 
Lemma~\ref{deg5} the coefficients of $f$ and $g$ satisfy the following
system of equations
\begin{eqnarray} \label{sys1}
2(a_2b_3-a_3b_2)+a_1b_4-a_4b_1&=&0  \nonumber \\
a_2b_4-a_4b_2&=&0 \\
a_3b_4-a_4b_3&=&0 \nonumber
\end{eqnarray}
Assume first that $a_4$ or $b_4 \neq 0.$
Without loss of generality we may assume that $a_4=1.$
Replacing $g$ by $g-b_4f$ we may assume $b_4=0.$
But then the system~(\ref{sys1}) implies $b_3=b_2=b_1=0$ and we obtain
\[det(t_0 q_1+t_1q_2)= \left( \begin{array}{ccc}
         t_0a_0+t_1b_0&t_0a_1&t_0a_2\\
         t_0a_1&t_0a_2&t_0a_3\\
         t_0a_2&t_0a_3&t_0
\end{array}\right)  \]
which has a double root $(0:1).$
Hence $a_4=b_4=0.$ Then~(\ref{sys1}) just says 
\begin{eqnarray} \label{sys2}
a_2b_3-a_3b_2&=&0
\end{eqnarray}
If $a_3=b_3=0,$
\[det(t_0 q_1+t_1q_2)= \left( \begin{array}{ccc}
         t_0a_0+t_1b_0&t_0a_1+t_1b_1&t_0a_2+t_1b_2\\
         t_0a_1+t_1b_1&t_0a_2+t_1b_2&0\\
         t_0a_2+t_1b_2&0&0
\end{array}\right)  \]
which has a triple root.
Finally, again without loss of generality, we may assume that
$a_3=1,b_3=0.$ Then~(\ref{sys2}) implies $b_2=0$ and
\[det(t_0 q_1+t_1q_2)= \left( \begin{array}{ccc}
         t_0a_0+t_1b_0&t_0a_1+t_1b_1&t_0a_2\\
         t_0a_1+t_1b_1&t_0a_2&t_0\\
         t_0a_2&t_0&0
\end{array}\right)  \]
which again has a double root $(0:1).$ This completes the proof of the theorem.\qed

\bigskip

Theorems~\ref{iso} and~\ref{sta} lead to the question whether the map $\Phi: M_D \longrightarrow \sex$ is
surjective. In order to analize this question, recall that we represented a
quartic divisor on the conic $s:  x_1^2-4x_0x_2$ by the binary quartic:
\[f(t_0,t_1)=a_0t_0^4+4a_1t_0^3t_1+6a_2t_0^2t_1^2+4a_3t_0t_1^3+a_4t_1^4 \] 
and the conic associated to it was given by the matrix
 $\left( \begin{array}{ccc}
         a_0&a_1&a_2\\
         a_1&a_2&a_3\\
         a_2&a_3&a_4
\end{array}\right) .$
We are identifying in this way the projective space $\P{4}=\P{4}(a_0:\dots:a_4)$ of
binary quartics with the $\p{4}$ of conics apolar to the conic $s.$
Similarly we consider $Gr(1,4)$
as a space of pencils of conics. Hence we get a map 
\[ Jac: Gr(1,4) \longrightarrow \P{6} \]
associating to every pencil of conics apolar to $s$ the Jacobian of its
associated pencil of binary quartics. According to the definition of the map 
$\Phi$ a binary sextic $j$ is in the image of $\Phi$  if and only if the
pre-image $Jac^{-1}(j)$(which consists of 5 pencils counted with
multiplicities) consists only of special pencils of conics, i.e. pencils
of conics $\lambda q_1+\mu q_2$ such that $det(\lambda q_1 + \mu q_2)$
admits a multiple root $(\lambda:\mu).$
Let $\mathcal {D} \subset Gr(1,4)$ denote the hypersurface of $Gr(1,4)$
given by the equation 
\[ discrim(det(\lambda q_1+\mu q_2))=0 \]
$Jac$ being a finite morphism, the image $Jac(\mathcal{D})$ is a hypersurface in $\P{6}.
$ If $U$ denotes the open set of $\P{6}$ parametrizing stable
sextics and $\rho:U \longrightarrow \sex $ the natural projection,
it is easy to see that also
\[ \mathcal{H}:= p(Jac(\mathcal{D}) \cap U)\]
is a hypersurface in $\sex.$ We obtain

\begin{prop}
The image $\Im m(\Phi)$ of $\Phi$ satisfies
\[ \sex \setminus \mathcal{H} \subset \Im m(\Phi) \subset \sex.	\]
 
\end{prop}

\bigskip

The following two examples show that both inclusions are strict, i.e.
$\Im m(\Phi)\neq \sex$ and $\sex \setminus \mathcal{H} \neq \Im m(\Phi).$
\begin{ex} \label{j2}
\rm{
Consider the binary sextic

\[ j_0=-16t_1^6-48 t_0^4t_1^2+128t_0^3t_1^3-48t_0^2t_1^4-16 t_1^6 \]

$Jac^{-1}(j_0)$ contains the following 3 pencils of conics 
$p_i=\lambda q_1^{i}+\mu q_2^i,$ $(i=1,2,3$) with

\[q_1^1= \left( \begin{array}{ccc}
         -1&0&0\\
         0&0&0\\
         0&0&1
\end{array}\right) \;\;\; q_2^1= \left( \begin{array}{ccc}
         -8&1&0\\
         1&0&1\\
         0&1&0
\end{array}\right)
\]

\[q_1^2= \left( \begin{array}{ccc}
         1&1&-1\\
         1&-1&0\\
         -1&0&1
\end{array}\right) \;\;\; q_2^2= \left( \begin{array}{ccc}
         0&1&0\\
         1&0&1\\
         0&1&0
\end{array}\right)
\]

\[q_1^3= \left( \begin{array}{ccc}
         -1&1&0\\
         1&0&0\\
         0&0&1
\end{array}\right) \;\;\; q_2^3= \left( \begin{array}{ccc}
         0&1&0\\
         1&0&1\\
         0&1&0
\end{array}\right)
\]
Since the Jacobian matrix of the map $Jac: Gr(1,4) \longrightarrow \P{6}$
at the points $p_1,p_2,p_3$ of $Gr(1,4)$ is of rank $5,6,5$ the pencils
$p_1$ and $p_3$ are in the ramification locus of $Jac.$ Hence 
$Jac^{-1}(j_0)
=\{p_1,p_2,p_3 \}.$ Now it is easy to check that $discr(det(\lambda q_1^i
+ \mu q_2^i))=0$ for $i=1,2,3.$ hence there is no Desargues
configuration $D_{\pi}$ with $\Phi(D_{\pi})=j_0.$ }
\end{ex}

\begin{ex} \label{j3}
\rm{
Consider the binary sextic

\[ j_1=48t_0^5t_1+48 t_0^2t_1^4 \]

$Jac^{-1}(j_1)$ contains the 2 pencils of conics 
$p_i=\lambda q_1^{i}+\mu q_2^i,$ $(i=1,2$) with

\[q_1^1= \left( \begin{array}{ccc}
         2 \sqrt[3]{-2}&1&\frac{1}{8}\sqrt[3]{4}\\
         1&\frac{1}{8}\sqrt[3]{4}&\frac{1}{4}\sqrt[3]{-2}\\
         \frac{1}{8}\sqrt[3]{4}&\frac{1}{4}\sqrt[3]{-2}&1
\end{array}\right) \;\;\; q_2^1= \left( \begin{array}{ccc}
         2 \sqrt[3]{-2}&-1&0\\
         -1&0&0\\
         0&0                                                                                        &0
\end{array}\right)
\]

\[q_1^2= \left( \begin{array}{ccc}
         1&0&1\\
         0&1&-\frac{1}{2}\\
         1&-\frac{1}{2}&0
\end{array}\right) \;\;\; q_2^2= \left( \begin{array}{ccc}
         0&0&1\\
         0&1&0\\
         1&0&0
\end{array}\right)
\]
One easily checks $discr(det(\lambda q_1^1+\mu q_2^1))=0$ and 
 $discr(det(\lambda q_1^2+\mu q_2^2))\neq 0.$
So $ \lambda q_1^2+\mu q_2^2$ comes from a Desargues configuration, whereas
 $\lambda q_1^1+\mu q_2^1$
does not.
This implies $j_1 \in \Im m(\Phi),$ but $j_1 \in \mathcal{H}.$}
\end{ex}
 
\section{Special Desargues configurations} 

Consider the 10 points in $\P{3}$

\[S_{12}:=\overline{e_1e_2} \cap\overline{e_3e_4e_5}=(1:1:0:0)\;\;\;\;\;\;
S_{34}:=\overline{e_3e_4} \cap\overline{e_1e_2e_5}=(0:0:1:1)\]
\[S_{13}:=\overline{e_1e_3} \cap\overline{e_2e_4e_5}=(1:0:1:0)\;\;\;\;\;\;
S_{15}:=\overline{e_1e_5} \cap\overline{e_2e_3e_4}=(0:1:1:1)\]
\[S_{14}:=\overline{e_1e_4} \cap\overline{e_2e_3e_5}=(1:0:0:1)\;\;\;\;\;\;
S_{25}:=\overline{e_2e_5} \cap\overline{e_1e_3e_4}=(1:0:1:1)\]
\[S_{23}:=\overline{e_2e_3} \cap\overline{e_1e_4e_5}=(0:1:1:0)\;\;\;\;\;\;
S_{35}:=\overline{e_3e_5} \cap\overline{e_1e_2e_4}=(1:1:0:1)\]
\[S_{24}:=\overline{e_2e_4} \cap\overline{e_1e_3e_5}=(0:1:0:1)\;\;\;\;\;\;
S_{34}:=\overline{e_3e_4} \cap\overline{e_1e_2e_5}=(1:1:1:0)\]
Recall that a Desargues configuration $D_{\pi}$ is special if and only if 
the plane $\pi$ contains one of the points $S_{ij}$ (and no point $e_i$).
The following theorem gives a characterization of special Desargues
configurations in terms of their associated binary sextics.
\bigskip

\begin{teo}
For a Desargues configuration $D_{\pi}$ the following statements 
are equivalent:\newline
(1) $D_{\pi}$ is a special Desargues configuration. \newline
(2) The binary sextic $\Phi(D_{\pi})$ admits a double point.
\end{teo} 

\textsc{Proof:} Suppose $\pi$ given by $f(z)= \sum_{i=1}^4 \alpha_i z_i=0$
is a nondegenerate plane in $\P{3},$ such that $\Phi(D_{\pi})$
admits a double point. According to Theorem~\ref{jac}
$\Phi(D_{\pi})$ is the sextic divisor on the von Staudt conic defined
by the complete intersection $I(f,g,h)$ with 
$g(z)= \sum_{i=1}^4 \alpha_i z_i^2=0$ and  $h(z)= \sum_{i=1}^4 \alpha_i z_i^3=0.$
A point on this complete intersection scheme is a double point if and only if
the rank of the Jacobian matrix $J(f,g,h)$ is $\leq 2$ at this point.
But

\[J(f,g,h)=
\left(
\begin{array}{c}
df\\
dg\\
dh  \end{array} \right)
=\left( \begin{array}{cccc}
         \alpha_1&\alpha_2&\alpha_3&\alpha_4\\
         2\alpha_1z_1&2\alpha_2 z_2&2 \alpha_3 z_3& 2 \alpha_4 z_4\\
         3\alpha_1z_1^2& 3 \alpha_2z_2^2& 3 \alpha_3 z_3^2&3 \alpha_4z_4^2
\end{array}\right)
\]
So $rk(J(f,g,h))
=rk \left( \begin{array}{cccc}
         1&1&1&1\\
         z_1& z_2& z_3& z_4\\
         z_1^2&z_2^2& z_3^2&z_4^2
\end{array}\right) \leq 2
$ if and only if 
\[det\{ \left( \begin{array}{cccc}
         1&1&1&1\\
         z_1& z_2& z_3& z_4\\
         z_1^2&z_2^2& z_3^2&z_4^2
\end{array}\right) \mbox{ with i}^{th}\mbox{ column omitted} \}=0\]
for $i=1,\dots,4.$ But these are Vandermonde determinants. Hence
$(z_1:\dots:z_4) \in \Phi(D_{\pi})$ is a double point if and only if it
satisfies the following system of equations

\[ \left\{ 
\begin{array}{ccc}
(z_1-z_2)(z_1-z_3)(z_2-z_3)&=&0\\
(z_1-z_2)(z_1-z_4)(z_2-z_4)&=&0\\        
(z_1-z_3)(z_1-z_4)(z_3-z_4)&=&0\\        
(z_2-z_3)(z_2-z_4)(z_3-z_4)&=&0          
\end{array}\right.
\]
The solutions of the system above are exactly the points in
$\P{3}$ with 3 equal coordinates or two pairs of equal coordinates.
So after normalizing one type of such points
is represented by $(1:1:1:\gamma)$ and the other by $(1:1:\gamma:\gamma)$
with $\gamma \in \complex.$ Such a point is contained in the complete intersection
$I(f,g,h)$ if and only if $\gamma=1$ or $\gamma=0.$
If $\gamma=1$ then $s_{\pi}$ contains the point $e_5,$ so $\pi$ is degenerate.
If $\gamma=0$ then the $(z_1: \dots :z_4)$ is one of the 10 points
 $S_{ij}$ as above. We conclude that $\Phi(D_{\pi})$ admits
a double point if and only if $D_{\pi}$ is special. \qed
\bigskip

\begin{rmk}
{\rm In the same way one can show that a Desargues configuration
$D_{\pi}$ admits 2 (respectively 3) lines containing 4 points if and only if 
$\Phi(D_{\pi})$ admits 2 (respectively 3) double points.}

\end{rmk}
\section{Degenerate Desargues configuration}
In this section we study an extension of 
 the map $\Phi: M_{D} \longrightarrow \sex$
to a holomorphic map $\overline{\Phi}:U
\longrightarrow \sex$ where $U$ is an open set
 with $M_{D} \subset U \subset \overline{M}_D,$
where the first inclusion is strict.
First we have to define the von Staudt conic 
for degenerate Desargues configurations.

\bigskip

Recall that the von Staudt conic of a non degenerate Desargues
configuration $D_{\pi}$ is defined by the fact  that the map associating 
to every point $p_{ij}$ in $D_{\pi}$ the line $l_{klm}$ with complementary
indices is a polarity.
This definition can be generalized verbatim to the case of generalized
Desargues configurations. One only has to add that if 2 points 
$p_{ij}$ and  $p_{kl}$ fall together, no line is associated to them.

To be more precise, let $D_{\pi}$ be of the first kind, i.e. $D_{\pi}$
contains exactly one of the points $e_i.$ For the sake of notational
simplicity we assume that $e_5 \in \pi$ (see Figure 3). In this case the map
associating to $p_{ij}$ $(i,j \neq 5)$
the line $\ell_{kl5}$ is just given by the involution $\imath$ on the pencil
of lines through $e_5$ given by the fact that $D_{\pi}$ is a limit of
Desargues configurations (see Section 1). This involution has 2 fixed lines
and their union is the von Staudt conic $s_{\pi}$.

Let now $D_{\pi}$ be of the second kind. Without loss of generality, we assume
the case of Figure 4, i.e. $e_4$ and $e_5 \in \pi.$ 
There are only 3 points apart from $e_4$ and $e_5$ namely $p_{12},p_{13}$ and
$p_{23}.$ To $p_{12}$ the line $\ell_{345}$ is associated, to $p_{13}$ the line
$\ell_{245}$ and to $p_{23}$ the line $\ell_{145}.$ But these 3 lines coincide
with $\overline{e_4e_5}.$
Hence the von Staudt conic $s_{\pi}$ in this case is the double line 
$\overline{e_4e_5}.$

Finally if $D_{\pi}$  is of the third kind, there are 3 points in the plane
$\pi$ to which no line is associated. This implies that the von Staudt conic
$s_{\pi}$ is the zero conic, i.e. the whole plane $\pi.$ Note that if $D_{\pi}$
is of the $i^{th}$ kind, the rank of the conic is $3-i.$ 
Analogously to Lemma~\ref{VS} we have:

\bigskip
\begin{lema}
Suppose $D_{\pi}$ is a degenerate Desargues configuration. The von Staudt 
conic $s_{\pi}$ of $D_{\pi}$ is given by the equation $\sum_{i=1}^4\alpha_iz_i^2=0 $
in the plane $\pi=\sum_{i=1}^4\alpha_iz_i=0.$
\end{lema}

\bigskip    
\textsc{Proof:}
One can prove this either in the same way as Lemma~\ref{VS}
using an analogous version of Reyes remark (see proof of Lemma~\ref{VS})
or just check it by computation:\newline
If $D_{\pi}$ is degenerate of the first kind with say $e_4 \in \pi,$ i.e.
$\pi$ is given by $\alpha_1z_1+\alpha_2z_2+\alpha_3z_3=0$ with
$\alpha_1,\alpha_2,\alpha_3,\sum_{i+1}^3 \alpha_i \neq 0,$
the two fixed lines of the involution on the pencil of lines with centre 
$e_4$ are

\begin{eqnarray}
&\ell_1:&(-\alpha_1\alpha_2\alpha_3-\alpha_2
\sqrt{D})z_2-(-\alpha_1\alpha_2\alpha_3+\alpha_3\sqrt{D})z_3=0\\ \nonumber
\mbox{ and }&&\\ \nonumber
&\ell_2:&(-\alpha_1\alpha_2\alpha_3+\alpha_2
\sqrt{D})z_2-(-\alpha_1\alpha_2\alpha_3-\alpha_3\sqrt{D})z_3=0 \nonumber
\end{eqnarray}
with $D=-\alpha_1\alpha_2\alpha_3(\alpha_1+\alpha_2+\alpha_3) \neq 0.$ One
only has to check that the ideal $(\sum_{i=1}^3 \alpha_i z_i,\ell_1.\ell_2)$
coincides with the ideal of $s_{\pi}.$

If $D_{\pi}$ is degenerate of the second kind with say $e_3,e_4 \in \pi$
one easily checks that $(\alpha_1 z_1+\alpha_2 z_2,
\alpha_1 z_1^2+\alpha_2 z_2^2)$
is the ideal of the double line $s_{\pi}=2\overline{e_3e_4}.$
Finally  if $D_{\pi}$ is degenerate of the third kind with say $e_1,e_2,e_3
\in \pi$ then $(\alpha_1z_1,\alpha_1z_1^2)$ is the ideal of the 
whole plane $\pi.$ \qed

\bigskip

In order to investigate a possible extension of the map $\Phi:M_D
\longrightarrow \sex,$ recall its definition. Consider a Desargues
configuration with von Staudt conic $s_{\pi}.$ The pencil of conics 
passing through a complete quadrangle of $D_{\pi}$ restricts to a pencil of
quartic divisors in $s_{\pi}$ whose Jacobian is a stable sextic $J$
independent of the choice of complete quadrangle and we defined
$\Phi(D_{\pi})=J.$ Note that $J$ was given by  $\sum_{i=1}^4\alpha_iz_i=
\sum_{i=1}^4\alpha_iz_i^2=\sum_{i=1}^4\alpha_iz_i^3=0,$
if $\pi:\sum_{i=1}^4\alpha_iz_i=0.$
Now if $D_{\pi}$ is degenerate of the first kind (respectively second)
kind, $D_{\pi}$ admits still 4 (respectively 3) complete quadrangles.
Restricting the pencil of conics passing through one of them to the von Staudt
conic $s_{\pi},$ its Jacobian is still a sextic divisor on $s_{\pi}.$  
The same proof as for Theorem~\ref{VS} also gives in this situation 
that the sextic divisor on $s_{\pi}$ is independent of the choice of the
complete quadrangle and given by the equations
\begin{eqnarray}\label{deg}
\sum_{i=1}^4\alpha_iz_i=
\sum_{i=1}^4\alpha_iz_i^2=\sum_{i=1}^4\alpha_iz_i^3=0,
\end{eqnarray}
if $\pi:\sum_{i=1}^4\alpha_iz_i=0.$
(Notice that the difference to the above situation is that one (or two)
of the constants $\alpha_1,\dots,\alpha_4,\alpha_1+\dots+\alpha_4$ 
vanishes.)
Now if $D_{\pi}$ is of the first kind, the von Staudt conic is the union of
two lines which intersect in a point, say $p_0 \in \pi$ and it is easy to
check that the divisor $j_{\pi}$ on $s_{\pi}$ given by 
equations~(\ref{deg})
is $6p_{0}.$ If $D_{\pi}$ is of the second kind, the von Staudt conic is a
double line $s_{\pi}=2 \ell_{\pi}$ and it is easy to see that the zero set 
of~(\ref{deg}) is the whole line $\ell_{\pi}.$
In the first case $j_{\pi}$ can be interpreted as the sextic with a 6-fold point
and in the second case not a sextic at all. In any case $j_{\pi}$ cannot (at
least in an obvious way) be interpreted as a semistable binary sextic.
This means that the geometric definition of the map $\Phi$ does not extend to
degenerate Desargues configurations.
\begin{rmk}
{\rm At first sight this seems to contradict the fact that the rational map 
$\Phi:\overline{M}_{D} \cdots \longrightarrow \sex$ is a morphism outside a
subvariety of codimension $\geq 2.$ The explanation comes from the fact that
the space of all binary sextics is not separated.
} 

\end{rmk}

\begin{center}
\vskip20pt
$\begin{array} {c}
\hbox{\footnotesize Departamento de Matem\'atica, UFMG}\\
\hbox{\footnotesize Belo Horizonte, MG 30161--970, Brasil.}
\\
\hbox{\footnotesize dan@mat.ufmg.br}
\end{array}$\quad\quad
$\begin{array} c

\hbox{\footnotesize Mathematisches Institut}\\
\hbox{\footnotesize Bismarckstr. $1\frac{1}{2},$ 91054, Erlangen}
\\

\hbox{\footnotesize lange@mi.uni-erlangen.de}
\end{array}$
\end{center}
\end{document}